\documentclass[reqno]{amsart}

\usepackage{hyperref}
\usepackage{subfigure, graphicx}

\newcommand{\pa}[2]{\frac{\partial #1}{\partial #2}}
\newcommand{\mr}[1]{\mathrm{#1}}
\newcommand{\mc}[1]{\mathcal{#1}}
\newcommand{\mf}[1]{\mathfrak{#1}}

\newcommand{\mbb}[1]{\mathbb{#1}}
\newcommand{\lr}{\leftrightarrow}
\newcommand{\xr}{\xrightarrow}

\DeclareMathOperator{\rank}{rank}
\DeclareMathOperator{\minor}{minor}
\DeclareMathOperator{\Ima}{Im}
\DeclareMathOperator{\Ker}{Ker}
\DeclareMathOperator{\Coker}{Coker}

\DeclareMathOperator{\Tr}{Tr}

\DeclareMathOperator{\Ad}{Ad}

\theoremstyle{plain}
\newtheorem{thm}{Theorem}[section]

\newtheorem{cnj}{Conjecture}

\theoremstyle{definition}
\newtheorem{dfn}{Definition}[section]
\newtheorem{exa}{Example}[section]

\theoremstyle{remark}
\newtheorem{rmk}{Remark}[section]

\begin{document}

\title{Euclidean geometric invariants of links in 3-sphere}

\author{Evgeniy V. Martyushev}

\address{South Ural State University, 76 Lenin avenue, 454080 Chelyabinsk, Russia}

\email{mev@susu.ac.ru}

\thanks{The work is partially supported by Russian Foundation for Basic Research, Grant no. 04-01-96010}

\begin{abstract}
We present a new link invariant which depends on a representation of the link group in $\mr{SO}(3)$. The computer calculations indicate that an
abelian version of this invariant is expressed in terms of the Alexander polynomial of the link. On the other hand, if we use non abelian
representation, we get the squared non abelian Reidemeister torsion (at least for some torus knots).
\end{abstract}

\maketitle

\section*{Introduction}\label{sec:intro}

In this paper we consider a new link invariant. Its construction is naturally divided into three main parts. First, on a given representation of the
link group we define a covering of 3-sphere branched along the link. Then, we map the covering space into 3-dimensional Euclidean space according to
the representation. In the last, algebraic part, we build an acyclic complex; the torsion of this complex is the main ingredient of our invariant.

Such an invariant was first constructed by I.G.~Korepanov for 3-manifolds in the paper~\cite{Kor01}. There the simplest version of the invariant was
considered corresponding to the trivial covering of a manifold. Calculations showed that this version raised to the $(- 1/6)$th power equals the
order of the torsion subgroup of the first homology group. However, when we use the universal cover, which corresponds to the unity of the
fundamental group, we get more interesting version of the invariant associated with the (abelian) Reidemeister torsion, see~\cite{KM02, Mar03}.

We have written a few computer programs which allow to calculate the invariant for a given manifold $M$ (or link $L$) and for a given representation
$\rho$ of the fundamental group (or the link group) in the group of orientation preserving motions of $\mbb{R}^3$. With the help of these programs we
made two conjectures. The first one states a connection of the abelian version of our invariant with the Alexander polynomial of link. The second
conjecture states that for a torus knot and non abelian representation of its group, our invariant is the squared non abelian Reidemeister torsion
investigated and calculated in~\cite{DubPhD}.

The paper is organized as follows. In section~\ref{sec:constr} we define the invariant for a given link and representation of its group. In
section~\ref{sec:calc} we propose an example of calculation of abelian and non abelian versions of the invariant for the trefoil knot. In the last
section we suggest our conjectures.

\textbf{Acknowledgements.} I am glad to thank I.G.~Korepanov for proposing me the problem and numerous helpful discussions and remarks.

\section{Constructing the invariant}\label{sec:constr}

\subsection{Simplicial moves}

\begin{figure}
 \centering
 \includegraphics[scale=0.25]{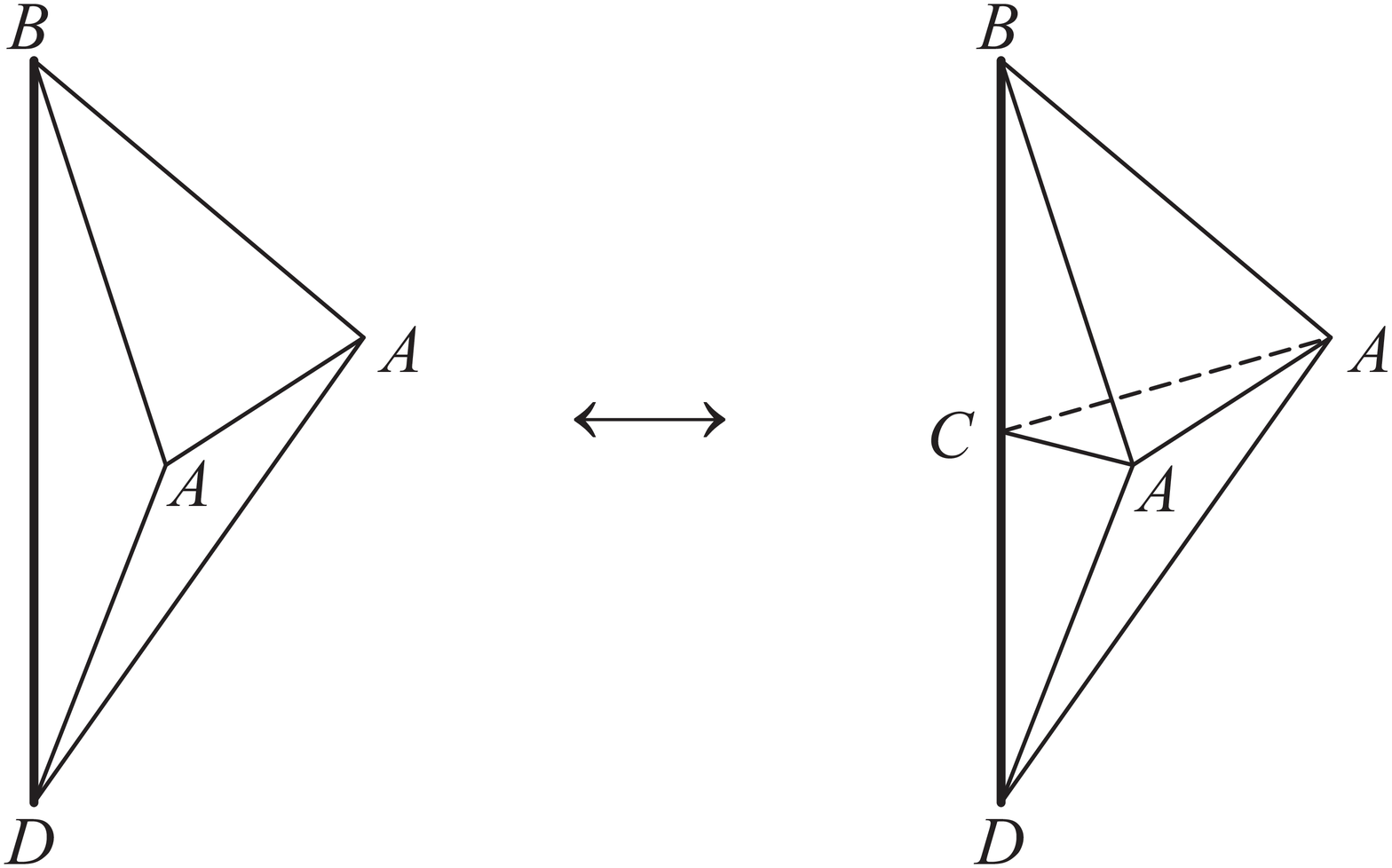}
 \caption{Moves $1 \lr 2$}\label{fig:move12}
\end{figure}

Let us describe a way of constructing the invariant following~\cite{Kor04} with minor changes in notation. Suppose we are given a link $L \subset
S^3$ endowed with a certain orientation. Let us consider a triangulation of $S^3$ satisfying the following conditions:
\begin{enumerate}
\item
the whole link $L$ lies on certain edges of the triangulation; \label{enu:first}
\item
for any tetrahedron in the triangulation, not more than two of its vertices belong to $L$;
\item
any edge $e$ of the triangulation either has two different vertices as its ends or, if its ends coincide, $e$ represents a meridian of the
corresponding link component. \label{enu:last}
\end{enumerate}
We need the edges with coinciding ends from condition~\ref{enu:last}) to define the following simplicial moves. Let an edge $BD$ lie on a certain
link component, and let there be a tetrahedron $BDAA$ in the triangulation, with its edge $AA$ representing a meridian of the corresponding link
component. The move $1 \to 2$ is defined as follows: take a point $C$ in the edge $BD$ and replace the tetrahedron $BDAA$ by two tetrahedra $BCAA$
and $CDAA$. The move $2 \to 1$ is the inverse to that (fig.~\ref{fig:move12}).

Using methods of the paper~\cite{Lick}, one can prove the following analogue of the Pachner theorem (see~\cite{Kor04} for details).
\begin{thm}\label{thm:Kor04}
A triangulation of sphere $S^3$ obeying the conditions \ref{enu:first}) -- \ref{enu:last}) can be transformed into any other triangulation obeying
the same conditions by a sequence of the following elementary moves:
\begin{itemize}
\item
Pachner moves $2 \lr 3$ and $1 \lr 4$. Such moves are not affect the edges lying on the link $L$ (however, the link may pass through edges and/or
vertices lying in the boundary of the transformed cluster of tetrahedra);
\item
moves $1 \lr 2$.
\end{itemize}
\end{thm}

\subsection{Triangulation of 3-sphere}

To define the invariant we need a triangulation of 3-sphere satisfying the conditions~\ref{enu:first}) -- \ref{enu:last}). We can construct such a
triangulation for arbitrary link as follows.

Strictly speaking, we are going to construct a \textit{pseudotriangulation} of 3-sphere. It differs from a triangulation in the proper sense in that
a simplex in a pseudotriangulation can appear several times in the boundary of a simplex of greater dimension. Clearly, any pseudotriangulation can
always be transformed into triangulation using a barycentric subdivision. However, in our case this would make no sense, since such a subdivision
significantly increases the number of vertices and edges and hence makes the calculations much more difficult.

So, let $r$ be the full amount of crossing points of the link diagram. Since the Euler characteristic of 2-sphere is 2, it follows that the link
diagram splits the plane by $r+2$ regions.

Denote by $C_i$ the $i$th crossing point of the link diagram. Let the plane with the link diagram be $z = 0$ in $\mbb{R}^3$. Take a point $O$ in a
region $z < 0$, and a point $P$ in $z > 0$. Consider a bipyramid over each of $r+2$ polygons with $O$ and $P$ as its apexes. Adding $r+2$ more edges
$OP$, which correspond to the $r+2$ bipyramids, and replacing $\mbb{R}^3$ with $S^3 \cong \mbb{R}^3 \cup \{\infty\}$, we obtain a triangulation of
$S^3$ consisting of $4r$ tetrahedra.

\begin{figure}
 \centering
 \subfigure[Initial triangulation]
 {\includegraphics[scale=0.15]{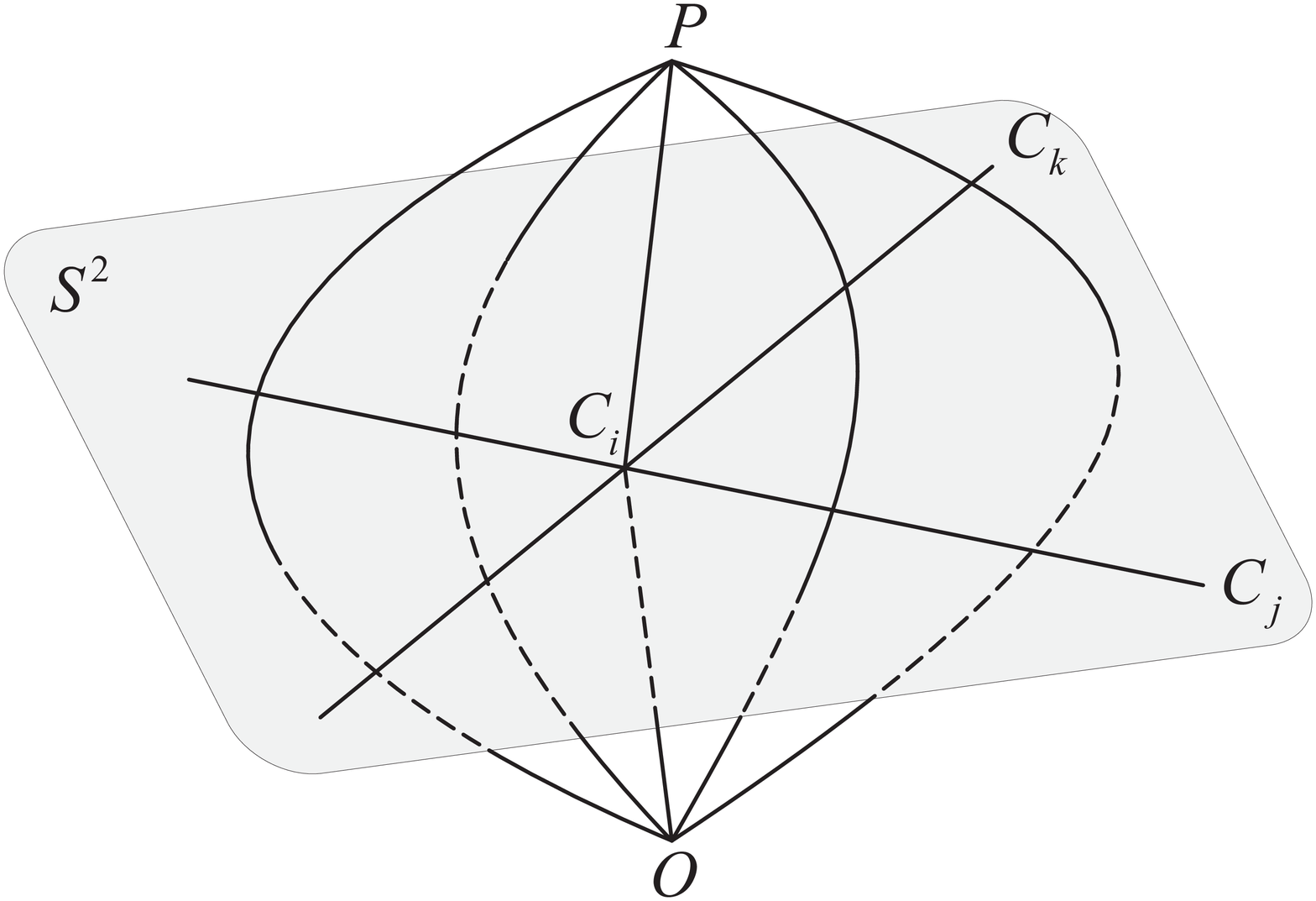}\label{fig:triangS3a}} \qquad
 \subfigure[Transformed triangulation (two edges $A_i O$ and two edges $B_i P$ are not depicted here)]
 {\includegraphics[scale=0.15]{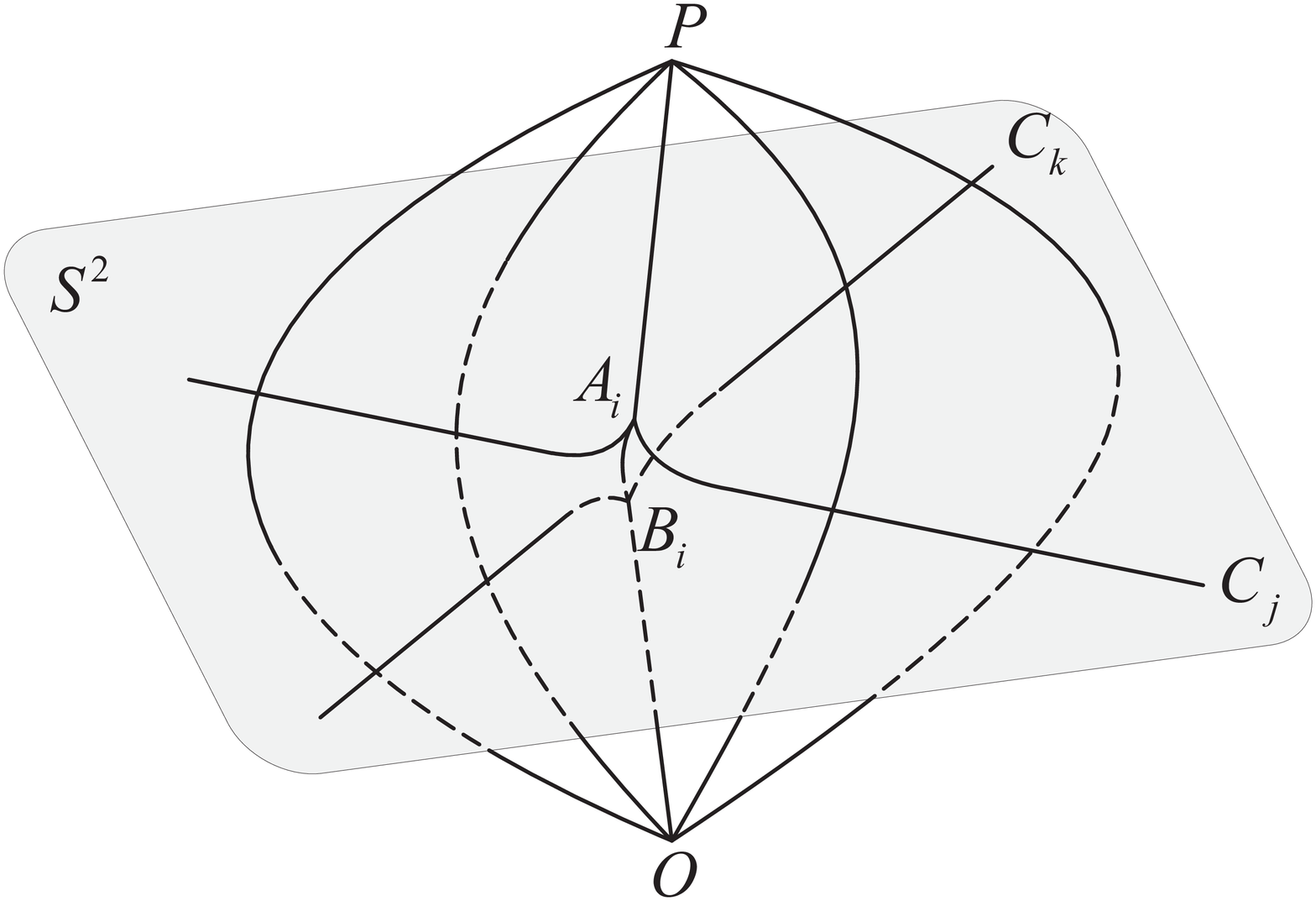}\label{fig:triangS3b}}
 \caption{Transforming 3-sphere triangulation in neighborhood of point $C_i$}
 \label{fig:triangS3}
\end{figure}

Further, we transform the obtained triangulation in a neighborhood of each point $C_i$ in such a way as depicted in figure~\ref{fig:triangS3}.
We suppose the point $A_i$ to be above $B_i$. Note that the tetrahedron $C_i C_j OP$ is replaced with $A_i C_j OP$, and the tetrahedron $C_i C_k OP$
with $B_i C_k OP$. Besides, four more tetrahedra $A_i B_i OP$ arise.

Doing this transformation for all $i = 1, \ldots, r$, we obtain a desired triangulation of $S^3$ with the link $L$ going through some of its edges.
Surely, the conditions~\ref{enu:first}) -- \ref{enu:last}) hold for this triangulation.

Thus, the set of vertices (0-simplices) of the triangulation consists of $2r + 2$ elements: $A_1, \ldots, A_r$, $B_1, \ldots, B_r$, $O$ and $P$.

The set of edges of the triangulation splits in four subsets:
\begin{enumerate}
\item $2 r$ edges lying on the link;
\item $r$ edges $A_i B_i$, where $i = 1,\ldots, r$;
\item $6 r$ edges $A_i O$, $A_i O$, $A_i P$, $B_i O$, $B_i P$, $B_i P$, where $i = 1,\ldots, r$;
\item $r + 2$ edges $O P$.
\end{enumerate}
Hence, the full amount of edges (1-simplices) in the triangulation equals $10r + 2$. It is easy to see that the full amount of 2- and 3-simplices is
$16r$ and $8r$ respectively.

\subsection{Coverings of 3-sphere branched along links}

\begin{dfn}\label{dfn:brcover}
Let $M$ and $N$ be a couple of triangulated closed orientable 3-manifolds and let a link $L$ be composed of some edges of manifold $N$.
\textit{Branched covering} along the link $L$ is such a continuous map $p \colon M \to N$ that the preimage $p^{-1}(L)$ is 1-dimensional subcomplex
in $M$ and the restriction $p$ to $M \setminus p^{-1}(L)$ is an ordinary covering.
\end{dfn}

Due to R.~Fox~\cite{Fox}, a branched covering is uniquely determined by its ordinary covering induced by restriction. Denote by $\rho \colon \pi L
\to E(3)$ the \textit{nontrivial} representation of the link group in the group of orientation preserving motions of 3-dimensional Euclidean space.
Then, for the normal subgroup $\Ker \rho$ there exists a unique covering $p_{\rho} \colon \tilde{S}^3 \to S^3$ branched along the link $L$ and its
multiplicity is $|\pi L / \Ker \rho|$.

\begin{rmk}\label{rmk:image}
Let $\pi L = \langle x_1, \ldots, x_r \mid R_1, \ldots, R_{r-1} \rangle$ be the Wirtinger presentation of the link group. A connected component of
the link diagram is called the \textit{overpass}. Recall that each overpass uniquely corresponds to a certain generator $x_i$. Therefore, if a vertex
$A$ belongs to the $i$th overpass, where $i = 1, \ldots, r$, then $x_i A = A$. It follows that the element $\rho(x_i)$ is represented by a rotation
(without translations) around some axis in $\mbb{R}^3$ in such a way that the conditions $\rho(R_1)$, \ldots, $\rho(R_{r-1})$ hold. The form of
relations $R_j$ implies that the rotation through an angle $\varphi_j$ around the axis occurs whenever the corresponding overpass belongs to the
$j$th link component.
\end{rmk}

\begin{exa}\label{exa:abel}
Let us consider a representation $\rho$ that sends each generator $x_i$ to the rotation through an angle $\varphi_j$ about the same fixed axis going
through the coordinate origin in $\mbb{R}^3$, provided that the generator $x_i$ goes around the $j$th link component. For this representation, $\Ker
\rho = [\pi L, \pi L]$ and $\Ima \rho$ is an abelian subgroup in $\mr{SO}(3)$ (such a representation is called \textit{abelian}). The corresponding
covering is called the \textit{universal abelian covering} branched along the link $L$.
\end{exa}

Let us fix the \textit{fundamental family} $\mc{F}$ for the covering space $\tilde{S}^3$, i.e. a family of simplices of $\tilde{S}^3$ such that over
each simplex of $S^3$ lies exactly one simplex of this family.

\subsection{Acyclic complex and link invariant}

We shortly remind basic definitions from the theory of algebraic complexes, see~\cite{Tur01} for details.

Let $C_0$, $C_1$, \ldots , $C_n$ be finite-dimensional $\mbb{R}$-vector spaces. We suppose that each $C_i$ is based, that is has distinguished basis.
Then, linear mapping $f_i \colon C_{i + 1} \to C_i$ can be identified with matrix.

\begin{dfn}\label{dfn:cmplx}
The sequence of vector spaces and linear mappings
\begin{equation}\label{eq:cmplx}
C = (0 \xr{} C_n \xr{f_{n-1}} C_{n-1} \xr{} \ldots \xr{} C_1 \xr{f_0} C_0 \xr{} 0)
\end{equation}
is called a \textit{complex} if $\Ima f_i \subset \Ker f_{i-1}$ for all $i = 1, \ldots , n - 1$. This condition is equivalent to $f_{i-1} f_i = 0$
for all $i$.
\end{dfn}

\begin{dfn}\label{dfn:homol}
The space $H_i(C) = \Ker f_{i-1} / \Ima f_i$ is called the \textit{$i$th homology} of the complex $C$.
\end{dfn}

\begin{dfn}\label{dfn:acycl}
The complex $C$ is said to be \textit{acyclic} if $H_i(C) = 0$ for all $i$. This condition is equivalent to $\rank f_{i-1} = \dim C_i - \rank f_i$
for all $i$.
\end{dfn}

Suppose that the sequence~\eqref{eq:cmplx} is an acyclic complex. Let $\mc{C}_i$ be an ordered set of basis vectors in $C_i$ and let $\mc{B}_i
\subset \mc{C}_i$ be a subset of basis vectors belonging to the space $\Ima f_i$.

Denote by ${}_{\mc{B}_i} f_i$ a nondegenerate transition matrix from the basis in space $\Coker f_{i+1} = C_{i+1}/\Ima f_{i+1}$ to the basis in space
$\Ima f_i$. By acyclicity, such a matrix really exists. Hence, ${}_{\mc{B}_i} f_i$ is a principal minor of the matrix $f_i$ obtained by striking out
the rows corresponding to vectors of $\mc{B}_{i+1}$ and the columns corresponding to vectors of $\mc{C}_i \setminus \mc{B}_i$.

\begin{dfn}\label{dfn:tors}
A quantity
\begin{equation}\label{eq:torsion}
\tau(C) = \prod\limits_{i = 0}^{n-1} (\det {}_{\mc{B}_i} f_i)^{(- 1)^{i+1}}
\end{equation}
is called the \textit{torsion} of acyclic complex $C$.
\end{dfn}

\begin{thm}[\cite{Tur01}]
Up to a sign, $\tau(C)$ does not depend on the choice of subsets $\mc{B}_i$.
\end{thm}

\begin{rmk}\label{rmk:changing}
The torsion $\tau(C)$ does depend on the distinguished basis of $C_i$. If one performs change-of-basis transformation in every space $C_i$ with
nondegenerate matrix $A_i$, then the torsion $\tau(C)$ is multiplied by
\[
\prod\limits_{i = 0}^n (\det A_i)^{(-1)^{i+1}}.
\]
\end{rmk}

\begin{dfn}\label{dfn:dual}
Given a finite-dimensional $\mbb{R}$-vector space $V$ with a distinguished basis $\mathbf{e}$, let $V^*$ be its \textit{dual} with a basis
$\mathbf{e}^*$. That is, if $e_j$ is a vector from $\mathbf{e}$ and $e^*_i$ is a vector from $\mathbf{e}^*$, then
\[
e^*_i (e_j) =
\begin{cases}
1, & i = j, \\
0, & i \neq j.
\end{cases}
\]
\end{dfn}

Let us return to construct the invariant. Given a branched covering $p_\rho$, we produce a continuous map $\Gamma \colon \tilde{S}^3 \to \mbb{R}^3$
as follows. Let $F$ be a vertex belonging to the triangulation of $S^3$. Then we place the orbit of $F$, i.e. the set $\{g \tilde{F} \mid g \in \pi
L\}$, in $\mbb{R}^3$ according to the following rule: if $\tilde{F}_2 = g \tilde{F}_1$, where $g \in \pi L$, then $\Gamma(\tilde{F}_2) = \rho(g)
\Gamma(\tilde{F}_1)$. Thus, we assign to each vertex of the triangulation the set of points in $\mbb{R}^3$ and these points are related by means of
motions from $E(3)$.

Further, each simplex of nonzero dimension is mapped into the convex linear shell of its vertex images in $\mbb{R}^3$. These shells may intersect
each other in arbitrary way but all the incidence relations are surely preserved.

After that, we can assign to each edge its Euclidean length and to each tetrahedron its volume. Besides, other Euclidean quantities (dihedral angles,
for instance) now make sense.

Note that all the tetrahedra entered in the triangulation of oriented space $\tilde{S}^3$ can be oriented \textit{consistently}, i.e., for every
tetrahedron, we can order its vertices up to even permutations. Under the mapping $\Gamma$, a tetrahedron either preserves this orientation or
changes it to opposite. In the first case we take the volume of the tetrahedron and all its dihedral angles with the sign $+$, in the second case
with the sign $-$.

\begin{rmk}
The vertices of fundamental family $\mc{F}$ are mapped in $\mbb{R}^3$ in arbitrary way. However, we require that the configuration of these vertices
in $\mbb{R}^3$ obeys the following conditions of general position:
\begin{itemize}
\item
the volumes of all tetrahedra from $\Gamma(\mc{F})$ are nonzero;
\item
the rank of $\left(\pa{\omega_i}{l_j}\right)$ possesses the maximal value at the point $\omega_i=0$, $\forall i$; here, $\omega_i$ is the defect
angle at $i$th edge, $l_j$ is the length of $j$th edge, $i$ and $j$ run over all the edges from $\mc{F}$.
\end{itemize}
\end{rmk}

Now we are going to define several vector spaces in order to combine them into the acyclic complex~\eqref{eq:acycl}.

\begin{dfn}
The set of elements $g \in E(3)$ such that $g h g^{-1} = h$ for all $h \in \Ima \rho$ is called the \textit{centralizer} of $\Ima \rho$.
\end{dfn}

It is a subgroup in $E(3)$ and we denote by $\mf{e}(3)_\rho$ its Lie algebra endowed with a distinguished basis. Thus,
\[
\mf{e}(3)_\rho = \{u \in \mf{e}(3) \mid \Ad_{\rho(h)} u = u, \forall h \in \pi L\}.
\]

Let $v \in \mc{F}$ be a vertex from the fundamental family. We define a vector space $(dx)_v$ at the point $\Gamma(v)$ as follows. Provided that
$\Gamma(v)$ does not belong to the link, $(dx)_v$ consists of column vectors $(dx_v\; dy_v\; dz_v)^T$, where $dx_v$, $dy_v$ and $dz_v$ are the
differentials of Cartesian coordinates of $\Gamma(v)$. Otherwise, $(dx)_v = u_x dx_v + u_y dy_v + u_z dz_v$, where $(u_x, u_y, u_z)$ is the unit
directing vector of the axis whereon the point $\Gamma(v)$ lies. Let us define $(dx) = \bigoplus\limits_{v \in \mc{F}} (dx)_v$. Here the direct
summation is over all the vertices of $\mc{F}$.

Suppose $e \in \mc{F}$ is an edge from the fundamental family. Let $L_e = \frac{1}{2} l_e^2$, where $l_e$ is the Euclidean length of edge $\Gamma(e)
\in \mbb{R}^3$. We define the vector space $(dL) = \bigoplus\limits_{e \in \mc{F}} (dL_e)$, where $(dL_e)$ is 1-dimensional space generated by
$dL_e$.

Let $\omega_e$ be the defect angle at $e$ and $\Omega_e = \frac{\omega_e}{l_e}$. Similarly, we define the space $(d\Omega) = \bigoplus\limits_{e \in
\mc{F}} (d\Omega)_e$.

We denote by $\mf{e}(3)^*_\rho$ and $(dx)^*$ the dual spaces to $\mf{e}(3)_\rho$ and $(dx)$ respectively (definition~\ref{dfn:dual}).

\begin{thm}[\cite{Kor02}, \cite{Kor04}]
The sequence
\begin{equation} \label{eq:acycl}
0 \xr{} \mf{e}(3)_\rho \xr{f_1} (dx) \xr{f_2} (dL) \xr{f_3 = f_3^T} (d\Omega) \xr{-f_2^T} (dx)^* \xr{f_1^T} \mf{e}(3)^*_\rho \xr{} 0
\end{equation}
is an acyclic complex. Here the superscript $T$ means, of course, matrix transposition.
\end{thm}

Since the bases in all the spaces of~\eqref{eq:acycl} are distinguished, the linear mappings $f_1$, $f_2$ and $f_3$ can be identified with matrices.

Denote by $\mc{C}_0$, $\mc{C}_1$ and $\mc{C}_2$ arbitrary ordered sets of basis vectors in the spaces $\mf{e}(3)_\rho$, $(dx)$ and $(dL)$
respectively. Let $\mc{B}_i$ be a subset of basis vectors belonging to the space $\Ima f_i$. Denote by ${}_{\mc{B}_i} f_i$ such a principal minor of
matrix $f_i$ that its rows correspond to the vectors from $\mc{C}_{i-1} \setminus \mc{B}_{i-1}$, and its columns correspond to the vectors from
$\mc{B}_i$. We suppose that $\mc{B}_3 = \mc{C}_2 \setminus \mc{B}_2$. Then, according to~\eqref{eq:torsion}, the torsion of complex~\eqref{eq:acycl}
is given by
\begin{equation}\label{eq:tors}
\tau = \frac{(\minor f_2)^2 (-1)^{\rank f_2}}{\minor f_3 \, (\minor f_1)^2}.
\end{equation}

Let $N$ be the amount of components of the link $L$. Put
\begin{equation}\label{eq:inv}
I_\rho(L) = \tau \cdot \frac{\prod\limits_{j = 1}^N (2 - 2\cos\varphi_j)^{n_j}}{\prod' l^2 \cdot \prod (-6V)}.
\end{equation}
Here $n_j$ is the amount of vertices in the fundamental family lying on $j$th component; $\varphi_j$ are the above mentioned angles (rotation around
the overpass belonging to $j$th component through the angle $\varphi_j$ determines a generator of the link group in the Wirtinger presentation). In
the denominator, the primed product is taken over whose edges from $\mc{F}$ which lie on the link, $V$ is a tetrahedron volume and the second product
is taken over all tetrahedra of $\mc{F}$.

\begin{thm}[\cite{Kor04}]
For a given representation $\rho$, the quantity $I_\rho(L)$ is an invariant of the link $L$. More precisely, it is independent of the choice of
$\mc{F}$ and $\Gamma$, and it does not change under the moves $2 \lr 3$, $1 \lr 4$ and $1 \lr 2$ (see theorem~\ref{thm:Kor04}).
\end{thm}

\section{Calculations for the trefoil knot}\label{sec:calc}

\begin{figure}
\centering
\includegraphics[scale=0.25]{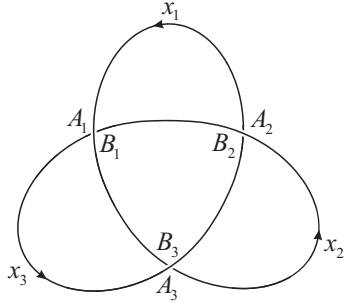}
\caption{Diagram of the trefoil}\label{fig:trefoil}
\end{figure}

\subsection{Case of abelian representation}

In figure~\ref{fig:trefoil} we depict an oriented diagram of the left-handed trefoil knot $3_1$. We are going to calculate the invariant
$I_\rho(3_1)$ for the abelian representation $\rho$ described in example~\ref{exa:abel}. There are three generators of the knot group denoted by
$x_1$, $x_2$ and $x_3$. They uniquely correspond to the three overpasses on the diagram and we denote these overpasses by the same letters.

Under $\Gamma$, every overpass $x_i$ is mapped into the line segment on the same axis $z$. A rotation around this axis through an angle $\varphi_1$
specifies an abelian representation of $\pi 3_1$ in $\mr{SO}(3)$. It follows that the Lie algebra $\mf{e}(3)_\rho$ is 2-dimensional, its basis
consists of infinitesimal rotation around the axis $z$ and infinitesimal translation along this axis. Therefore, $\mc{C}_0 = \mc{B}_0 = \{d\varphi_z,
dz\}$.

We place the vertices $A_i$ and $B_i$ of the fundamental family in $(0, 0, i+2)$ and $(0, 0, i-1)$ respectively, vertices $O$ and $P$ in $(1, 0, 0)$
and $(0, 1, 0)$ respectively.

The set $\mc{C}_1$ (set of basis vectors in $(dx)$) include the differentials of Cartesian coordinates of the vertices $O$, $P$ and the differentials
of $z$'s coordinate of the vertices $A_i$ and $B_i$ for $i = 1, 2, 3$. Choose
\[
\mc{B}_1 = \{dz_{A_1}, dy_O\}.
\]
Then,
\[
\det {}_{\mc{B}_1} f_1 = \begin{vmatrix} 0 & 1 \\ x_O & 0
\end{vmatrix} = - 1.
\]

The dimension of space $(dL)$ is 32: 10 basis vectors correspond to the matrix $f_2$, and 22 vectors correspond to the matrix $f_3$. Choose
\[
\begin{split}
\mc{B}_2 = \{dL_{B_3, A_1}, dL_{A_1, x_2B_2}, dL_{B_2, A_3}, dL_{A_3, x_1B_1}, dL_{B_1, A_2},
\\ dL_{A_1, x_3O}, dL_{A_1, P}, dL_{B_1, O}, dL_{x_1B_1, P}, dL_{O, P}\},
\end{split}
\]
where, for example, $x_2B_2$ means the action of element $x_2 \in \pi 3_1$ on the vertex $B_2$.

Computing the determinants of ${}_{\mc{B}_2} f_2$ of size $10 \times 10$ and ${}_{\mc{B}_3} f_3$ of size $22 \times 22$, we find
\[
\begin{split}
\det {}_{\mc{B}_2} f_2 &= 1440, \\
\det {}_{\mc{B}_3} f_3 &= - \frac{(1 - \cos \varphi_1)^4 \, (1 - 2\cos \varphi_1)^4}{167961600 \, \cos^{12} \varphi_1 \, (2\cos^2 \varphi_1 - 1)^6}.
\end{split}
\]

Then, substituting the found values for principal minors in~\eqref{eq:tors}, we obtain the torsion and then, by formula~\eqref{eq:inv}, the
invariant:
\begin{equation}\label{eq:IabelT23}
I_\rho(3_1) = - \frac{(2 - 2 \cos \varphi_1)^2}{(1 - 2 \cos \varphi_1)^4}.
\end{equation}

\subsection{Case of non abelian representation}

Let the trefoil group $\pi 3_1$ be presented by two generators $a$, $b$ and the only relation $a^2 = b^3$. We have
\begin{equation}\label{eq:x_i}
x_1 = a b^{-1}, \qquad x_2 = b^{-1} a, \qquad x_3 = a^{-1} b a b a^{-1},
\end{equation}
where, as above, $x_1$, $x_2$ and $x_3$ are the Wirtinger generators and also corresponding overpasses (fig.~\ref{fig:trefoil}).

Let us consider a representation $\rho \colon \pi 3_1 \to \mr{SO}(3)$ such that $\Tr(\rho(a)) = 1 + 2\cos(\pi)$, $\Tr(\rho(b)) = 1 +
2\cos(\frac{2\pi}{3})$. In other words, we fix two axes going through the coordinate origin in $\mbb{R}^3$. Then, $\rho(a)$ is a rotation around one
of these axes through the angle $\pi$, and $\rho(b)$ is a rotation around another axis through the angle $\frac{2\pi}{3}$. We denote by $\phi$ the
angle between these axes.

Using formulas~\eqref{eq:x_i}, one can see that under $\Gamma$ the overpasses $x_1$, $x_2$ and $x_3$ are mapped into the line segments on three axes
going through the coordinate origin. The angles between these axes are the same in pairs. The corresponding generators of the trefoil group are
represented by rotations around these axes through an angle $\varphi_1$. Simple calculation shows that this angle is associated with $\phi$ by the
relation:
\[
\cos \varphi_1 =\frac{3}{2} \cos^2 \phi - 1.
\]

Further, for the described representation we use the same triangulation as for the abelian case.

The centralizer of the subgroup $\Ima \rho$ is trivial now and hence, $\mf{e}(3)_\rho = 0$. It follows that
\[
\mc{B}_0 = \mc{B}_1 = \emptyset.
\]
Set
\[
\begin{split}
\mc{B}_2 = \{dL_{B_3, A_1}, dL_{A_1, x_2B_2}, dL_{B_2, A_3}, dL_{A_3, x_1B_1}, dL_{B_1, A_2}, dL_{A_1, B_1}, \\
dL_{A_1, x_3O}, dL_{A_1, O}, dL_{A_1, P}, dL_{B_1, O}, dL_{x_1B_1, P}, dL_{B_2, P}\}.
\end{split}
\]
Recall that $x_2B_2$ means an action of an element $x_2 \in \pi 3_1$ on the vertex $B_2$.

Further, we calculate the determinants of ${}_{\mc{B}_2} f_2$ of size $12 \times 12$ and ${}_{\mc{B}_3} f_3$ of size $20 \times 20$. Finally, with
the help of~\eqref{eq:tors} and~\eqref{eq:inv}, we get
\begin{equation}\label{eq:InonabelT23}
I_\rho(3_1) = 4.
\end{equation}

\section*{Discussion}

We see that the invariant $I_\rho(3_1)$ for the abelian representation $\rho$ is equal to
\[
- \left|\frac{e^{i\varphi_1} - 1}{\Delta_{3_1}(e^{i\varphi_1})}\right|^4,
\]
where $\Delta_{3_1}(t) = t^2-t+1$ is the Alexander polynomial for the trefoil knot. Calculations for other knots and links give similar results.
Therefore, we propose the following
\begin{cnj}\label{cnj:links}
Let $L \subset S^3$ be a link with $N$ components, $\Delta_L(t_1, \ldots, t_N)$ its Alexander polynomial, and $\rho \colon \pi L \to \mr{SO}(3)$ an
abelian representation of the link group. Then,
\begin{equation*}
I_\rho(L) =
\begin{cases}
- |\Delta_L(e^{i\varphi_1})|^{- 4} \cdot
(2 - 2 \cos \varphi_1)^2, & N = 1 \\
- |\Delta_L(e^{i\varphi_1}, \ldots, e^{i\varphi_N})|^{- 4}, & N > 1. \\
\end{cases}
\end{equation*}
\end{cnj}

\bigskip

Note that the trefoil knot is a particular case of so-called torus knot.
\begin{dfn}
Let $\mu$ and $\lambda$ be a meridian and parallel of standard embedded 2-torus and let $p$ and $q$ be a couple of coprime integers. Then, a closed
curve $p\mu + q\lambda$ is called the \textit{torus knot} of type $(p, q)$ and is denoted by $T(p, q)$.
\end{dfn}

Recall that the group of torus knot $T(p, q)$ has a presentation $\pi T(p, q) = \langle a, b \mid a^p = b^q \rangle$. Consider a representation
$\rho_{j, k} \colon \pi T(p, q) \to \mr{SO}(3)$ such that $\rho_{j, k}(a)$ is a rotation around a certain axis $x$ through the angle $\frac{2\pi
j}{p}$ and $\rho_{j, k}(b)$ is a rotation around another axis $x'$ through the angle $\frac{2\pi k}{q}$, where $1\leq j\leq \lfloor p/2\rfloor$ and
$1\leq k\leq \lfloor q/2\rfloor$. The angle between $x$ and $x'$ is nonzero.

\begin{cnj}\label{cnj:ITpq}
The invariant $I_{j, k}(T(p, q))$ for the representation $\rho_{j, k}$ looks like:
\begin{equation}
I_{j, k}(T(p, q)) = \frac{1}{(pq)^2} \left(4 \sin \frac{j\pi}{p} \sin \frac{k\pi}{q}\right)^4.
\end{equation}
\end{cnj}
Note that this conjecture agrees with~\eqref{eq:InonabelT23} for the trefoil $T(2, 3)$. Note also that the non abelian Reidemeister torsion twisted
by the action of group $\mr{SU(2)}$ for the $T(p, q)$ is equal to (see~\cite[p. 113]{DubPhD})
\[
-\frac{1}{pq} \left(4 \sin \frac{j\pi}{p} \sin \frac{k\pi}{q}\right)^2
\]
(if we use the meridian curve to compute the torsion, which gives a natural basis for the twisted $H^1$; the author would like to thank J.~Dubois for
pointing out this fact). Thus, our non abelian invariant seems to be the square of the non abelian Reidemeister torsion.

\bibliographystyle{amsplain}

\begin{thebibliography}{99}

\bibitem{DubPhD}
Dubois J. Torsion de Reidemeister non abelienne et forme volume sur l'espace des repr\'{e}sentations du groupe d'un n{\oe}ud: Ph.D. thesis /
Universit\'{e} Blaise Pascal, 2003.

\bibitem{Fox}
Fox R.H. Covering spaces with singularities // Algebraic Geometry and Topology: A Symposium in Honor of S.Lefschetz. Princeton Math. Series, 1957.
Vol.~12. P.~243--257.

\bibitem{Kor01}
Korepanov I.G. Invariants of PL manifolds from metrized simplicial complexes // J.~Nonlin. Math. Phys., 2001. Vol.~8. P.~196--210.

\bibitem{Kor02}
Korepanov I.G. Euclidean 4-simplices and invariants of four-dimensional manifolds: II. An algebraic complex and moves $2 \leftrightarrow 4$ // Theor.
Math. Phys., 2002. Vol.~133. P.~1338--1347.

\bibitem{Kor04}
Korepanov I.G. Euclidean tetrahedra and knot invariants // Proceedings of the Chelyabinsk Scientific Center, 2004. Vol.~24. P.~1--5.

\bibitem{KM02}
Korepanov I.G., Martyushev~E.V. Distinguishing three-dimensional lens spaces $L(7,1)$ and $L(7,2)$ by means of classical pentagon equation //
J.~Nonlin. Math. Phys., 2002. Vol.~9. P.~86--98.

\bibitem{Lick}
Lickorish W.B.R. Simplicial moves on complexes and manifolds // Geometry and Topology Monographs, 1989. Vol.~2. P.~299--320.

\bibitem{Mar03}
Martuyshev E.V. Euclidean simplices and invariants of three-manifolds: a modification of the invariant for lens spaces // Proceedings of the
Chelyabinsk Scientific Center, 2003. Vol.~19. P.~1--5.

\bibitem{Tur01}
Turaev V.G. Introduction to combinatorial torsions. Boston:~Birkhauser, 2000.~--~144p.

\end{thebibliography}

\end{document}